\documentclass[12pt]{article}

\input epsf
\usepackage{graphicx}
\usepackage{amssymb}
\usepackage{epstopdf}

\newtheorem{theorem}{Theorem}[section]
\newtheorem{lemma}[theorem]{Lemma}
\newtheorem{proposition}[theorem]{Proposition}

\newtheorem{defi}[theorem]{Definition}

\title{ Infinite sequence of fixed point free pseudo-Anosov homeomorphisms.}
\author{J\'er\^ome Los, \\
Universit\'e de Provence, LATP, UMR CNRS 6632.}

\begin{document}

\maketitle

{\bf Abstract.} We construct an infinite sequence of pseudo-Anosov homeomorphisms without fixed points and leaving invariant a sequence of orientable measured foliations on the same topological surface and the same stratum of the space of Abelian differentials. The existence of such a sequence shows that all pseudo-Anosov homeomorphisms fixing orientable measured foliations cannot be obtained by the Rauzy-Veech induction strategy.\\

{\em 2000 Mathematics Subject Classification.} {\footnotesize Primary: 37E30.  Secondary: 32G15, 57R30, 37B10}\\
{\em Key words and phrases.} {\footnotesize Pseudo Anosov homeomorphisms, Rauzy-Veech induction, fixed points, train-track maps}

\section{Introduction}

This work started after a discussion between  A.Avila and P.Hubert. The question was :
"Is there a pseudo-Anosov homeomorphism fixing an orientable measured foliation and without fixed points of negative index?"

It turns out that one such example already existed in the literature due to P.Arnoux and J.C.Yoccoz, [A.Y], [A]. It was constructed for a very different purpose and enables to build some special exemples on other surfaces but no general technics was available to build families of examples and in particular in the same stratum of the space of measured foliations.

Existence or non existence of fixed points is an interesting question in it's own right but why this particular setting?

One motivation comes from interval exchange transformations (IET) that naturally define topological surfaces $S$ together with orientable measured foliations. It also defines a natural transformation within the set of IET, the so-called {\em Rauzy-Veech induction} that has been widely studied over the years, in particular because it gives a natural relationship between the combinatorics of IET and some orbits of the Teichm\"uller flow (see for instance [Ra], [Ve]). 
Some of the IET are called {\em self-similar} and correspond to periodic loops of the Rauzy-Veech induction or to periodic orbits of the Teichm\"uller flow in moduli space.
They define elements in the mapping class group $\mathcal{M} (S)$ of the surface that are very often pseudo-Anosov, according to Thurston's classification theorem (see [Th], [FLP] for definitions). 
All these pseudo-Anosov homeomorphisms obtained from the Rauzy-Veech induction strategy share one property, in addition to fixing orientable measured foliations: they all have a fixed point and a fixed separatrix starting at this fixed point. In other words they all have a fixed point with negative index. The above question is thus rephrased as: Is there some pseudo-Anosov homeomorphisms with orientable invariant measured foliations that do not arise from the Rauzy-Veech induction strategy?
Our first result gives a positive answer to that question. This was not a surprise because of Arnoux-Yoccoz's example. 

The interest is more on the general construction, based on a tool called the {\em train-track automata} that has not been defined in full generality in the literature even though it exists in several forms and since a long time, see for instance, [Mo86], [PaPe] and more recently in [KLS].
These automata allow to construct, in principle, all pseudo-Anosov homeomorphisms fixing a measured foliation on a given stratum and are closely related to the {\em train-track complex} that U.Hamenst\"adt [Ha] has recently studied. These train-track automata will not be needed here, we will instead use a more elementary notion of {\em splitting sequence morphisms} that are obtained from the $R,L$-splitting sequences described by Papadopoulos and Penner in [PaPe] by an algebraic adaptation. In [PaPe] the aim was to describe measured foliations combinatorially.
Here it is a tool for explicitly constructing elements in the mapping class group by composition of elementary train-track morphisms. Observe that the Rauzy-Veech induction can be interpreted as a special class of splitting sequences on a special class of train-tracks.
Constructing just one example of pseudo-Anosov homeomorphisms without fixed point, using  splitting sequences, was a surprise and it became evident that such examples should be very rare on a given stratum. The natural conjecture was that at most finitely many such examples could exist on a given stratum and yet we obtain :

\begin{theorem}
There exists infinite sequences of pseudo-Anosov homeomorphisms with orientable invariant measured foliations on the same stratum of the space of measured foliations and without fixed points of negative index.

\end{theorem}

The existence of an infinite sequence of pseudo-Anosov homeomorphisms without fixed points of negative index is counter intuitive at several level. The general belief, prior to these examples, was that on a given stratum, at most, finitely many such pseudo-Anosov homeomorphisms could exist. This would have been a justification for the Rauzy-Veech induction strategy to be used for counting periodic orbits of the Teichm\"uller flow in moduli space. This counting problem is very delicate and several asymptotic results has been obtained recently, in particular by Eskin-Mirzakhani [EM] and Bufetov [Bu].
The examples obtained in this paper indicate that another strategy should be developed for a more precise counting. 
The non finiteness is also counter intuitive at the growth rate level. Indeed for a sequence of pseudo-Anosov homeomorphisms whose growth rate goes to infinity it was expected that the number of fixed points should grow, these examples show this is not the case.

The paper is organised as follows. First the basic properties of train-tracks and train-track maps are reviewed in section 2, in particular with respect to the properties that we want to control: the fixed points and the orientability. Then the notion of splitting sequence morphisms will be discussed in section 3. Finally the particular sequence leading to the examples will be described and analysed.

It's a great pleasure to thanks Pascal Hubert and Erwan Lanneau for their support and questions during this work. I would like to thanks Lee Mosher who made a constructive comment on a previous version of this paper.

\section{Train track maps and dynamics.}

In this section we review some basic properties of train-track maps that represents pseudo-Anosov homeomorphisms. In particular the relationship between the fixed points of the pseudo-Anosov homeomorphisms on the surface and those of the train-track maps. For more details about train-tracks, laminations and foliations the reader is referred to [PH] and [FLP].

\subsection{Train-track maps.}

A train-track on a surface $S$ is a pair $(\tau, h)$, where :\\
$\bullet$   $\tau$ is a graph, given by it's collection of edges $E(\tau)$ and vertices $V(\tau)$, together with a {\em smooth structure}, that is a partition at each vertex $v\in V(\tau)$, of the set of incident edges $St(v) = In(v) \bigsqcup Out(v)$, each subset, called a {\em side}, being non empty (see Figure 1).

\begin{figure}[htbp]
 \centering
\includegraphics[height=40mm]{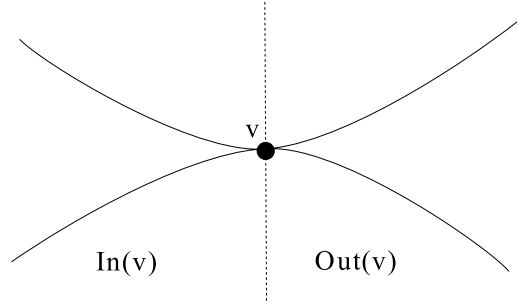}
\caption{ The smooth structure at a vertex.}

\label{fig:1}
\end{figure}

$\bullet$  $h : \tau \rightarrow S$ is an embedding that preserves the smooth structure and such that $S - h(\tau)$ is a finite union of discs with more than 3 {\em  cusps} on its boundary and annuli (if $S$ has boundaries and/or punctures) with more than one cusp on one boundary, the other boundary (or puncture) being a boundary component (or a puncture) of the surface $S$. A {\em cusp} here is given by a pair of adjacent edges, with respect to the embedding $h$, in the same side $In(v)$ or $Out(v)$ at a vertex.

$\bullet$  A regular neighborhood $N(h(\tau))$ on $S$ is a subsurface (with boundaries) and 
$h$ defines an embedding
$\hat{h} :  \tau \rightarrow N(h(\tau))$. The neighborhood $N(h(\tau))$ admits a retraction \\
$\rho :  N(h(\tau))  \rightarrow \tau$ that is a homotopy inverse of $\hat{h}$ and 
$\{ \rho^{-1} (t); t \in \tau \}$ is a foliation of $N(h(\tau))$ called the 
{\em tie foliation}. A leaf of the tie foliation that maps under $\rho$ to a vertex is called a 
{\em singular leaf}.

\begin{figure}[htbp]
\centerline {\includegraphics[angle=0,height=60mm]{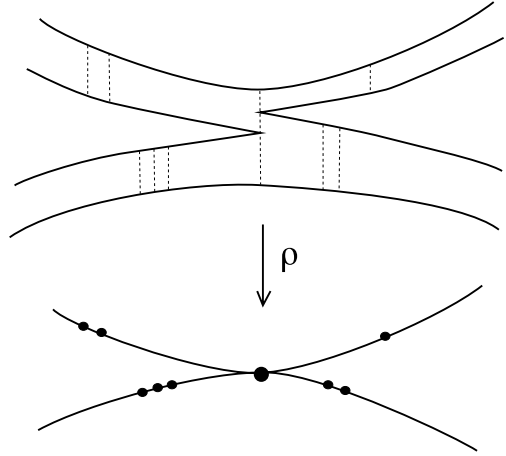} }
\label{fig:2}
\caption{ A regular neighborhood and the tie foliation.}
\end{figure}

$\bullet$  A map $\varphi : \tau \rightarrow  \tau '$ between two train-tracks $(\tau, h)$ and  $(\tau ', h ')$  on $S$ that satisfies :\\
(i)  $\varphi$ is cellular, i.e. maps vertices to vertices and edges to edge paths,\\
(ii) $\varphi$ {\em preserves} the smooth structure,
 is called a {\em train-track morphism}. It is a {\em train-track map} if $\tau$ and $\tau '$ are isomorphic (as train-tracks).\\
The last property (ii) means the following:\\
 - If $\varphi (v) = w$, $  (v,w) \in V(\tau)$ then the edges issued from one side ($In(v)$ or $Out(v)$) are mapped to edge paths starting at $w$  by edges on one side.\\
- If an edge path $\varphi (e)$ {\em crosses} a vertex $v$ then it enters the neighborhood of $v$ by one side and exits by the other side.

$\bullet$   If $f : S \rightarrow S$ is a homeomorphism, then a train-track map 
$\varphi : \tau \rightarrow  \tau$ is a {\em representative} of $[f] \in \mathcal{M}(S) $ if the following diagram commutes, up to isotopy:\\
\centerline{$\begin{array}{ccc}
\tau & \stackrel{h}{\longrightarrow} & S \\
{  {\varphi } }{\hspace{1mm}}{\downarrow} &                          & {\downarrow}{\hspace{1mm}}{ {f} }\\
\tau & \stackrel{h}{\longrightarrow} & S
\end{array}$}
In addition, the image graph $f \circ h (\tau)$ is supposed to be embedded into $N(h(\tau))$ and transverse to the tie foliation.\\
Such a triple $(\tau , h, \varphi)$ is called a {\em train-track representative} of the mapping class element 
$[f] \in \mathcal{M}(S)$.\\
For train-track maps there is a natural {\em incidence matrix} 
$M(\tau, \varphi)$ whose entries are labelled by the edge set $E (\tau)$ and 
$M(\tau, \varphi)_{(e, e')} $ is the number of occurrences of $e'^{\pm 1}$ in the edge path 
$\varphi(e)$. This definition requires the edges to be oriented but the orientation is arbitrary. This matrix does not depends upon the embedding $h$.\\
The following proposition is due to Thurston and the proof to Papadopoulos and Penner [PaPe].

\begin{proposition}
Any pseudo-Anosov homeomorphism admits a train-track representative whose incidence matrix is irreducible and non periodic. $\square$
\end{proposition} 

There is a converse direction result that requires an additional property to be checked, called the {\em "singularity type condition"}. The first formulation of this condition was given in the non published preprint version of the book [CB] and several equivalent conditions are described in the non published monograph by L.Mosher [Mo03]. The Casson-Bleiler version has been used in several published papers, for instance in [Ba].
The condition can be described as follows:
The {\em singularity type} of a measured foliation $ (\mathcal{F}, \mu)$ with $m$ singularities $\{s_1, ..., s_m\}$ is a m-tuple of integers $ [k_1, ..., k_m]_{\mathcal{F}}$, where $k_i$ is the number of separatrices at $s_i$.\\
A train-track $(\tau, h)$ on $S$ has {\em singularity type}  $ [k'_1, ..., k'_r]_{\tau}$ if the components of $S-h(\tau)$ are $r$ discs $\{\Delta_1, ... \Delta_r\}$ and $\Delta_j$ has $k'_j$ cusps on its boundary. This is the definition for closed surfaces, the adaptation for punctured surfaces is obvious.
If a measured foliation $ (\mathcal{F}, \mu)$ is {\em carried} by a train-track $(\tau, h)$ 
(see [PH] or [Mo03]) then we say that $ (\mathcal{F}, \mu)$ and $(\tau, h)$ have 
{\em identical singularity type} if $ [k_1, ..., k_m]_{\mathcal{F}} =  [k'_1, ..., k'_r]_{\tau} $ 
as non ordered ($m = r$)-tuples.
A reformulation of the Casson-Bleiler result is the following:

\begin{proposition}
If a train-track representative $(\tau, h, \varphi)$  of $[f] \in \mathcal{M}(S)$ has an irreducible and non periodic incidence matrix $ M(\tau, \varphi)$ and if the measured foliation 
$ (\mathcal{F}, \mu)_{(\tau, h, \varphi)}$ that is canonically obtained from $(\tau, h, \varphi)$ have identical singularity type, then $[f]$ is pseudo-Anosov.$\square$
\end{proposition} 

This result will be used in the last section to check our explicit examples. For completeness let us describe briefly the construction of the measured foliation $ (\mathcal{F}, \mu)_{(\tau, h, \varphi)}$ and one method that enables to check the singularity type condition.\\
A train-track representative $(\tau, h , \varphi)$ with an irreducible incidence matrix defines a canonical rectangle partition and a measured foliation of the subsurface $N(h( \tau ) )$ via the classical "highway construction" of Thurston (see [PH] for details).  The set of transverse measures (widths) on the rectangles is given by the positive eigenvector of the largest eigenvalue of the transpose matrix ${} ^{t}M(\tau, \varphi)$.
For simplicity we assume that the ambient surface $S$ is closed. The above construction defines a measured foliation of the subsurface $N (h(\tau))$. In order to extend this measured foliation to $S$ we need to apply a "zipping" construction that is essentially the "Veech zippered rectangles" operation where the zipping parameters (lengths) are computed directly from the map $\varphi$. 
In this closed surface case all the components of 
$S - N(h( \tau ) )$ are discs $\Delta_i$ bounded by curves $\partial_i$ with $k'_i$ cusps.
These boundary curves are represented by cyclic edge paths in $\tau$. Since 
$(\tau, h , \varphi)$ is a train-track representative of a surface homeomorphism then 
$\varphi_{| \partial_{i}}$ is a cyclic edge path that is homotopic to some boundary curve
$\partial_j$. The boundary components $\partial_i$ are thus homotopicaly permuted under $\varphi$ and the map $\varphi_{| \partial_{i}}$ has some periodic points along some sides of the cusped closed curve 
$\partial_i$, called {\em boundary periodic points}. We will see in our explicit examples of section 4 that these periodic points are easy to find.
A simple test for the singularity type condition is given by :

\begin{proposition}
 If the map $\varphi$ of the train-track representative $(\tau, h , \varphi)$  has exactly one boundary periodic point on each side of the boundary curves $\partial_i$ then the singularity type condition is satisfied by $(\tau, h , \varphi)$.
 \end{proposition}
 
 {\em Proof.}
 This {\em single boundary periodic point} condition 
implies that on each $\partial_i$ with $k'_i$ cusps there are exactly $k'_i$ boundary periodic points. The unique path that connects two consecutive boundary periodic points along $\partial_i$
has to cross a single cusp.  We observe that it is an {\em irreducible periodic Nielsen path} (see [J] and [BH]). Each boundary curve $\partial_i$ is thus a concatenation of $k'_i$ periodic irreducible Nielsen paths.
The zipping operation mentioned above is now easy to describe, this is just a (sequence of) {\em folding operations} (see [BH]), from the cusp to the boundary periodic points along a periodic Nielsen path.
In this case with only one boundary periodic point on each side of $\partial_i$ the zipping operation defines a measured foliation $(\mathcal{F}, \mu)_{(\tau, h , \varphi)}$ that is canonically defined from  $(\tau, h , \varphi)$. For this measured foliation the singularity type is naturally identified with the one of $(\tau, h)$ and the singularity type condition is satisfied. In particular the beginnings of the separatrices at the singularities are in bijection with the irreducible Nielsen paths. In addition the permutation of these separatrices under the pseudo-Anosov homeomorphism is given by the computable permutation of the Nielsen paths under $\varphi$ .
$\square$

\subsection{Dynamics of train-track maps, orientability.}

The relationship between the dynamics of a train-track map representing a pseudo-Anosov mapping class and the dynamics of the pseudo-Anosov homeomorphism is rather well understood (see for instance [Bo] or [Lo]). In this paragraph we only review the fixed point question that is much easier.

\begin{lemma}
A train-track map with irreducible incident matrix has a fixed point on an edge 
$e\in E(\tau)$ if and only if the edge path $\varphi (e)$ contains  
$e^{\pm 1}$ ( i.e. $M(\tau, \varphi)_{(e, e)} \neq 0 $).
\end{lemma}

{\em Proof.}
 The first observation is that condition (i) of a train-track map implies the existence of a Markov partition for $\varphi$. The classical combinatorial dynamics (see [ALM] for example) for maps with a Markov partition implies immediately the result. $\square$
 
The next step is to relate the fixed points of the train-track map representing a pseudo-Anosov mapping class with the fixed points of the pseudo-Anosov homeomorphism.

\begin{lemma}
If $f : S \rightarrow S$ is a pseudo-Anosov homeomorphism that is represented by a train-track map
 $(\tau, h, \varphi)$ with irreducible incidence matrix then for each fixed point $x$ of $f$ that is not a singular point of the invariant foliation there is a unique fixed point $y$ of $\varphi$  so that $\rho (x) = y$.
The possible missing fixed points of $f$ with respect to $\varphi$ are in one to one correspondance with the boundary components of $N(h(\tau))$ bounding discs on $S$, that are invariant up to homotopy under
 $\varphi$. In addition, each such boundary component is a union of periodic Nielsen paths between boundary periodic points and the dynamics of these Nielsen paths (under $\varphi$) reflects the dynamics of $f$ restricted to the separatrices at the fixed singular points.
\end{lemma}

{\em Proof.} A sketch of proof for the first statement goes as follow. The construction, described above, of the measured foliation $(\mathcal{F}, \mu)_{(\tau, h , \varphi)}$ from $(\tau, h, \varphi)$ gives in fact a Markov partition for the pseudo-Anosov homeomorphism $f$. If a fixed point of $f$ is not a singularity then it belongs to the interior of one rectangle of the constructed partition. The retraction $\rho$ then maps the stable segment that contains the fixed point, within the corresponding rectangle, to a point on an edge of $\tau$. This point is fixed for $\varphi$ by the commutative diagram of the definition. The proof of the last statement relating boundary periodic points and periodic Nielsen paths for $\varphi$ with singular fixed points for $f$ is in fact contained in the proof of Proposition 2.3, more details can be found for instance in [FL].$\square$

\begin{figure}[htbp]
\centerline {\includegraphics[angle=0,height=35mm]{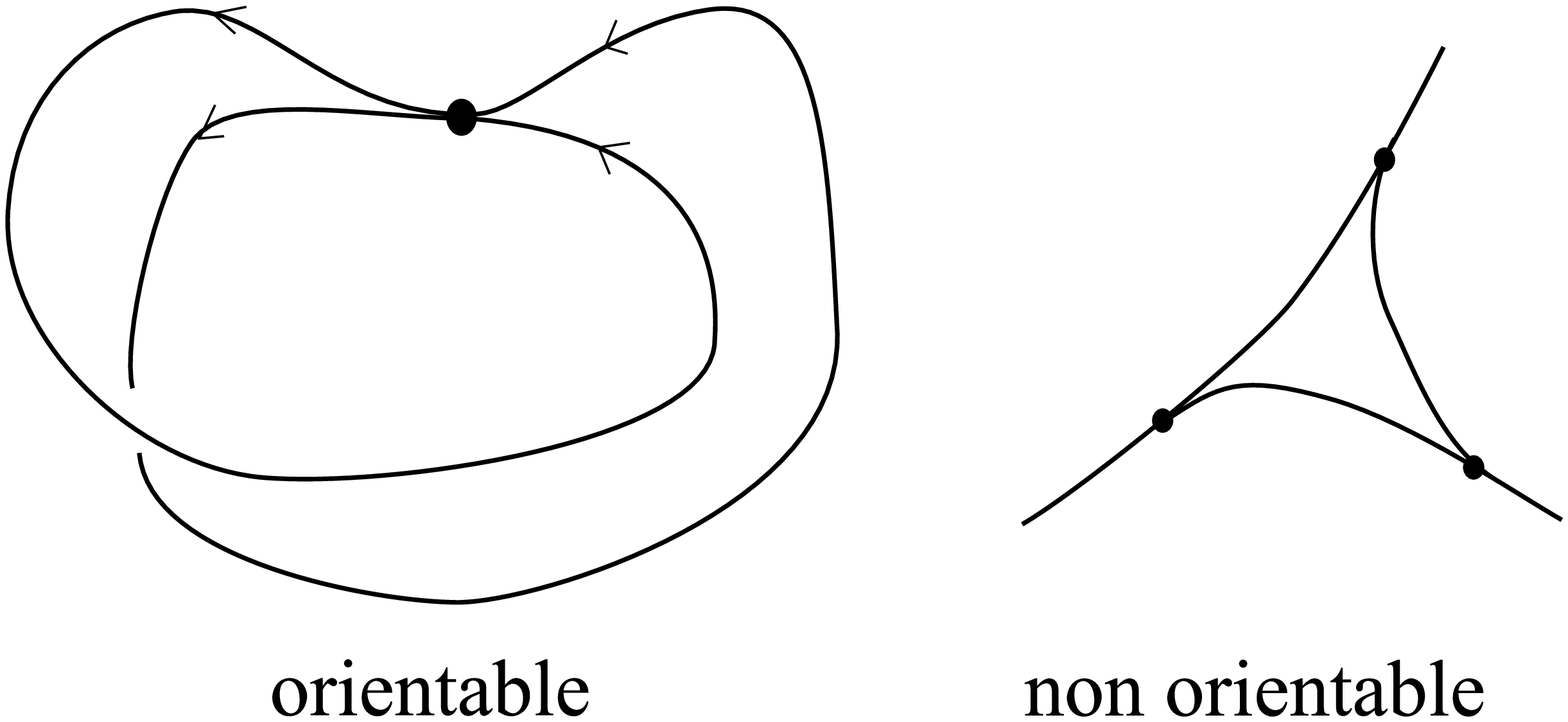} }
\label{fig:3}
\caption{ Orientable and non orientable train tracks.}
\end{figure}

The notion of an orientable foliation or measured foliation is classical.
We want to see this property at the train track level. The notion of being 
{\em carried} by a train track for a measured foliation (or a measured lamination) is also classical (see [PH]). If a measured foliation $(\mathcal{F}, \mu)$ is carried by $\tau$ then each leaf $\Lambda$ of  
$\mathcal{F}$ projects, under the retraction $\rho$, to an edge path that is {\em compatible} with the smooth structure in the sense of  property (ii) of a train-track map. If $(\mathcal{F}, \mu)$ is orientable then each leaf has a natural orientation and defines an orientation for the edge path $\rho(\Lambda)$. Each oriented edge path $\rho(\Lambda)$ are compatible with each others and induces an orientation for the edges of $\tau$. These orientations are compatible with the smooth structure, i.e. at each vertex the edges are oriented from the In side to the Out side (or the opposite). If a train-track admits a global orientation that satisfies this local compatibility property at each vertex it is called {\em orientable} and we have the obvious:

\begin{lemma}
If a train-track $(\tau, h)$ on $S$ is orientable then every measured foliation carried by $(\tau, h)$ is orientable.$\square$
\end{lemma}

\section{Splitting morphism sequences.}

In this section we introduce a general construction that produces train-track maps and the corresponding mapping class elements fixing measured foliations in a given stratum of the space of measured foliations.
The splitting operation is as old as the measured train-tracks, going back to the original paper by Williams [Wi]. The first result that relates pseudo-Anosov homeomorphisms or more precisely measured foliations invariant under a pseudo-Anosov homeomorphism to combinatorial splitting sequences is probably due to Papadopoulos and Penner [PaPe]. Splitting sequences have been used recently by Hamenst\"adt [Ha] as a main tool to study for instance the curve complex. In this section we introduce the notion of {\em splitting morphisms}.

Let $\varphi : \tau \rightarrow \tau$ be a train-track map representing 
$f : S \rightarrow  S$ then $h\circ \varphi : \tau \rightarrow N(h(\tau))$
is an embedding that is isotopic to $f\circ h$ and we assume that 
$f\circ h (\tau)$ is embedded in $N(h(\tau))$ transversely to the tie foliation.
Let us consider a small regular neighborhood 
$N' ( f\circ h (\tau) )$, small enough to be embedded in $N(h(\tau))$. This neighborhood is homeomorphic to $N(h(\tau))$ by definition and thus $N(h(\tau)) - N' ( f\circ h (\tau) ) $ is a union of annuli with the same number of cusps on each boundary component.
For each cusp $C'_i$ of $N' ( f\circ h (\tau) )$ there exists a cusp $C_i$ of $N(h(\tau))$ and a path 
$\gamma_i$ in $N(h(\tau)) - N' ( f\circ h (\tau) ) $ that connects $C'_i$ to $C_i$ and is transverse to the tie foliation. Such a path is called a {\em splitting path}. The following lemma first appeared in [PaPe]. 
\begin{lemma}
Every train-track map representing an element $[f]$ in the mapping class group is obtained by a finite sequence of "elementary" splitting operations followed by a relabeling homeomorphism.
\end{lemma}

\begin{figure}[htbp]
\centerline {\includegraphics[angle=0,height=40mm]{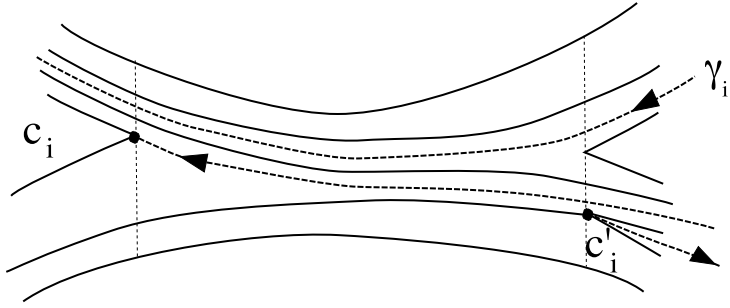} }
\label{fig:4}
\caption{ Splitting paths.}
\end{figure}

{\em Proof.} If we cut $N(h(\tau))$ along the finite collection of finite splitting paths $\gamma_i$ we obtain a new surface $N_1(\tau, h, \varphi)$ that is obviously homeomorphic with $N(h(\tau))$. This observation implies the result. $\square$

Our goal is now to make the notion of {\em elementary splitting} operation explicit, as well as the relabeling homeomorphism. Then we shall use these splitting operations for the reverse purpose, namely for constructing train-track maps and the corresponding mapping class elements.\\
In their paper Papadopoulos and Penner [PaPe]  were interested to splitting sequences at "generic" train- tracks, i.e. with only valency 3 vertices. For these sequences a very simple coding with 2 symbols is available, the so-called R,L sequences. More recently Hamenst\"adt added the constraints that $(\tau, h)$ represents foliations in the {\em principal} stratum, i.e. with only triangles for the complementary components.
We will not use these restrictions here. In particular for orientable foliations we need to allow foliations to belong to more degenerate strata. We will also encourage the train-tracks to have higher valency vertices, in order to simplify the definition of the splitting operations. We want the elementary operations to be train-track morphisms, given by cellular maps between two different train-tracks that we call {\em splitting morphisms}, to this end we need the train-tracks to be {\em labelled}.

Observe that splitting operations appeared for interval exchange transformations and the Rauzy-Veech induction strategy is nothing but a very special splitting sequence at a very special class of train-tracks.

A {\em labeling} is a map $\epsilon : E(\tau) \rightarrow \mathbb{A}$, where 
$\mathbb{A}$ is a finite alphabet. The map $\epsilon$ just names each edge.
A splitting operation is defined for a given $N(h(\tau))$ at a specific cusp $C_i$. We want to fixe the alphabet once and for all, and in particular that the number of edges is constant. We also want to avoid changing from one stratum to another. To this end we impose the following constraints:

(a) $\tau$ has at least a vertex of valency larger than 3.\\
(b) A splitting operation is not applied at a vertex of valency 3.\\
(c) At a vertex with more than 2 edges on one of it's sides, only  the splittings at the 
{\em extreme} cusps and through the {\em extreme} edges are allowed.\\
(d) A splitting path do not connect a cusp of $N(h(\tau))$ to another cusp of  $N(h(\tau))$.\\

The notion of an extreme cusp or edge is clear from the context. 
A splitting operation is a cutting along a path $\gamma$ that connects a cusp $C_i$ of $N(h(\tau))$ to a point in the interior of a singular leaf of the tie foliation (condition (d)). A splitting is {\em elementary} if the path $\gamma$ do not intersects singular leaves except at it's extreme points. The splitting operation then defines a new neighborhood $N(h'(\tau '))$ that is embedded into $N(h(\tau))$ and thus, via the retraction $\rho$, a map 
$ S_{C_i} : (\tau, h) \leftarrow (\tau' , h' )$ (it will be written 
$\tau  \leftarrow \tau'  $ if no confusion is possible).
Under the above conditions the following properties are obvious:

(i) Every train-track satisfying (a) admits some splitting operations.\\
(ii) The structure of the complementary components is preserved by (d).\\
(iii) The number of vertices and edges is preserved by (b), (c) and (d).\\
(iv) The map $ S_{C_i}$ is injective on vertices by (b) and (c).\\
(v) After a splitting operation the new train-track  $\tau'$ satisfies (a).\\

Since the number of edges is fixed under a splitting operation we fixe once and for all the alphabet 
$\mathbb{A}$.\\ 
Let us consider a labelled train-track $(\tau, h, \epsilon)$. We want to define explicitly the splitting morphism $ S_{C_i}$, as a cellular map using the labels. To this end we give an arbitrary orientation to the edges in order to describe the images $ S_{C_i}(e)$, for every edge $e\in E(\tau)$, as an edge path in 
$\tau$. A particular splitting morphism at a cusp $C_i$ is a local operation that changes the train-track in a neighborhood of the vertex $v$ where $C_i$ is based. If $C_i$ is on one side of $v$ there are, at most, two possible splitting morphisms at $C_i$, by condition (c).
Let us consider the particular case where $v$ has valency 4 with two edges on each side (see Figure  5), the general case is an obvious adaptation. Figure 5 shows a particular labeling and how the labeling is transformed under an elementary splitting. As observed above the transformation of the graph takes place in a neighborhood of $v$.\\
($*$) The natural convention for the labeling of $\tau'$ is to keep the same names as in 
$\tau$ for all the edges that still exist after the operation. Thus only one edge could be renamed after the operation but the alphabet is fixed and then only one name is available.\\
 This convention explains the unambiguous labeling of Figure 5. Observe that this convention also induces an orientation for the edges of the new train-track if the initial edges were oriented. In particular if $\tau$ is orientable then the new train-track $\tau'$ is also orientable and each edge has a well defined orientation. This convention is valid and unambiguous for any splitting morphism obtained via the conditions (a)-(d).

\begin{figure}[htbp]
\centerline {\includegraphics[angle=0,height=55mm]{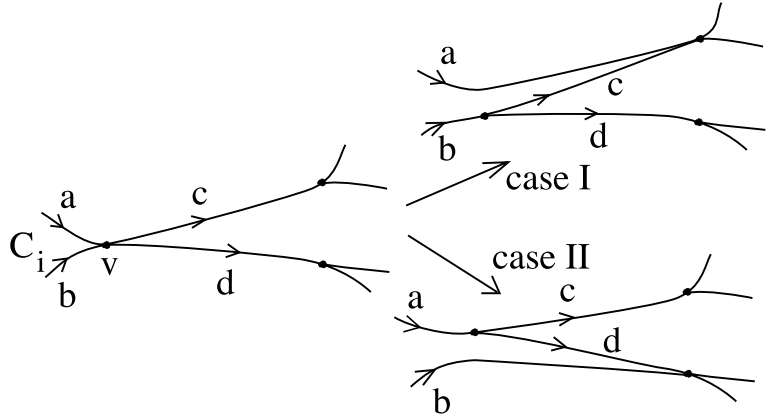} }
\label{fig:5}
\caption{  The splitting morphisms.}
\end{figure}

In the particular cases of Figure 5, i.e. with two edges (labelled c or d) on the opposite side of $C_i$, there are two possible splittings (case I and II in Figure 5) and the two morphisms are given by the following maps :\\
\centerline{
$S_{C_i} : \tau ' \rightarrow \tau$}\\
 \centerline{ $ a  \rightarrow  a.c,  \hspace{3mm}   x  \rightarrow x, $ for all $x \neq a$, in case I}\\
\centerline{ $ b  \rightarrow  b.d,  \hspace{3mm}  x  \rightarrow x, $ for all $x \neq b$, in case II.}

Lemma 3.1 can be rephrased as: 

Any train-track map $\varphi :\tau \rightarrow \tau$ representing an element in the mapping class group is obtained as a composition of splitting morphisms 
$\{  S_{C_{i_1}}, \cdots ,   S_{C_{i_n}} \} $, followed by a train-track automorphism  
$\alpha$:\\
\centerline{ $ \tau  \stackrel{ S_{C_{i_1}} }{\leftarrow} \tau_{1}
 \stackrel{ S_{C_{i_2}} }{\leftarrow} \tau_{2} \cdots 
 \stackrel{ S_{C_{i_n} }}{\leftarrow} \tau_{n}\simeq \tau
\stackrel{\alpha}{\leftarrow} \tau $. }

Along the sequence of splitting morphisms, each $\tau_j$ has a well defined labeling, obtained inductively from the one of $\tau$ by the convention ($*$). In particular $\tau_n$
 is homeomorphic with $\tau$ but the labeling is a priori different. The train-track automorphism $\alpha$ induces, in particular, a relabeling.
A given train-track map admits usually many splitting morphisms decompositions and this is a serious practical problem.

Let us now make an obvious but crucial observation for the fixed point problem.
The first observation is that the train-track automorphism might not be unique, this is the case if $\tau$ admits a symmetry, i.e. a non trivial automorphism group.\\
The splitting morphisms that are given above, together with Lemma 2.4, shows that if $\tau'$ is homeomorphic with $\tau$ and if the relabeling is the trivial one then many fixed points exist. Therefore to avoid fixed points for the train-track map it is highly desirable for the train-track $\tau$ to admit a non trivial automorphism group.

\section{Construction of examples.}

\subsection{A preliminary example.}

From the previous sections we are looking for a train-track map 
$\varphi : \tau \rightarrow \tau$ satisfying the following properties :\\
(i) $\tau$ is orientable,\\
(ii) $\tau$ has a non trivial automorphism group,\\
(iii) $\varphi$ has no fixed points.

Of course the first two conditions are quite easy to satisfy and the real difficulty is to find a map $\varphi$, as a sequence of splitting morphisms, satisfying (iii).
The first example is obtained by an experiment on the labelled train-track of Figure 6.

\begin{figure}[htbp]
\centerline {\includegraphics[angle=0,height=75mm]{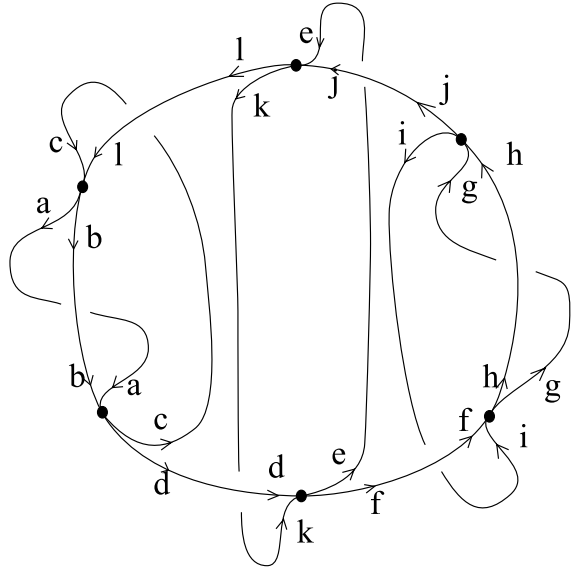}}
\label{fig:6}
\caption{ The initial labelled train-track.}
\end{figure}

On this train-track the orientability is obvious to check, at each vertex there is an orientation of the incident edges that is compatible with the In/Out partition of the smooth structure. The symmetry is of order two which means that the graph has two parts that are exchanged by the symmetry (an involution). With the labeling of  Figure 6 this involution is  given by the map :
$ \{ a, b, c, d, k, l\}\overrightarrow{\leftarrow} \{ g, h, i, j, e, f \}$, this notation means that the edge labelled $\{a\}$ is mapped to the edge labelled $\{g\}$ (and$\{g\}$ to $\{a\}$) and so on.
 Any splitting sequence produces a sequence of maps that increases the length of the image of some edge. Our main strategy is to avoid any image to be "too long" on it's symmetric part.

The difficulty for constructing examples is easily expressed by observing that the train-track 
$\tau$ of Figure 6 has 12 cusps and the maximal number of possible splittings at any train-track along the sequence is 24. Thus the number of possibilities at each step is of order 24 and the growth is exponential with respect to the length of the sequence.

For describing the sequence we introduce a simple notation. First notice that starting from a labelled and oriented train-track as in Figure 6 implies that after each splitting step, the new train-track has a canonical orientation and labeling by (*). A given edge $x\in E (\tau )$, being oriented, has an {\em initial} vertex 
$i(x)$ and a {\em terminal} vertex $t(x)$. The splitting morphisms of Figure 5 can be described as "sliding" the terminal part of $"a"$ {\em over} the initial part of $"c"$
in case I or the terminal part of $"b"$ over the initial part of $"d"$ in case II. For this reason the two possible morphisms are denoted 
$ \frac{t (a)}{i(c)}$ or $ \frac{t (b)}{i(d)}$, this notation is unambiguous for a labelled and oriented train track and gives at the same time the cusp and the morphism.

The very first sequence we want to describe on the labelled train-track of Figure 6, with this notation, is the following one:
\begin{center}
$  \mathcal{S}_1  =   \Big\{ \frac{i (b)}{t(l)} ; \frac{t (b)}{i(d)} ; \frac{t (k)}{i(f)} ; \frac{i (k)}{t(j)} ;
\frac{t (e)}{i(l)} ; \frac{t (c)}{i(a)} ; \frac{i (c)}{t(a)} ; \frac{i (e)}{t(d)} ;
\frac{i (h)}{t(f)} ;  \frac{t (h)}{i(j)} ; \frac{t (f)}{i(g)} ; \frac{i (j)}{t(g)} \Big\}$.
\end{center}

The surprise is that such a short sequence gives back to the same topological train-track. The composition of all these elementary splittings gives the image train-track as shown by Figure 7.
The figure should be understood as embedded in the regular neighborhood of the topological  (i.e. without the labels) train-track $\tau$ of Figure 6.

We observe that this image train-track is indeed the same topological train-track but with a different labeling, given by Figure 8, and of course a different embedding. In addition, since the automorphism group of $\tau$ has order two, there are two ways to relabel the train-track by a train track automorphism, this is shown by the two "identifications" in Figure 8.

\begin{figure}[htbp]
\centerline {\includegraphics[angle=0,height=90mm]{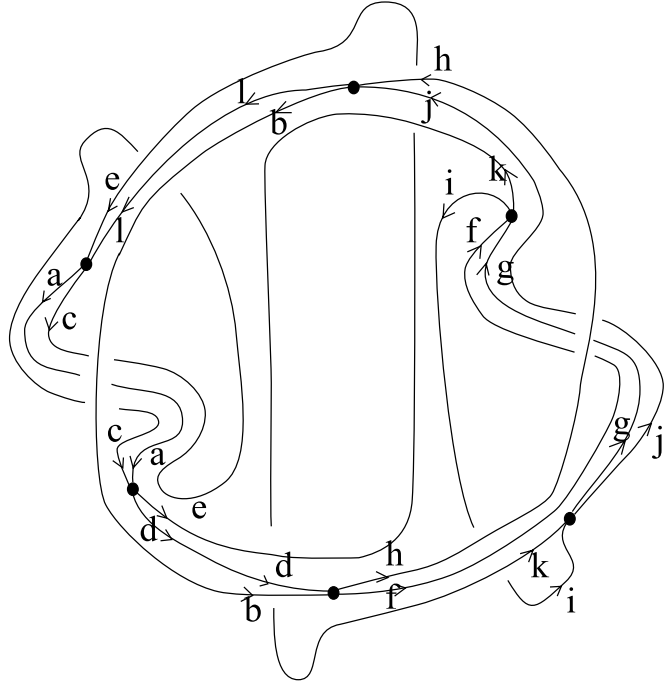} }

\label{fig:7}
\caption{ The image train track, after the sequence $  \mathcal{S}_1$ of splittings, in the regular neighborhood of $\tau$.}
\end{figure}

\begin{figure}[htbp]
\centerline {\includegraphics[angle=0,height=70mm]{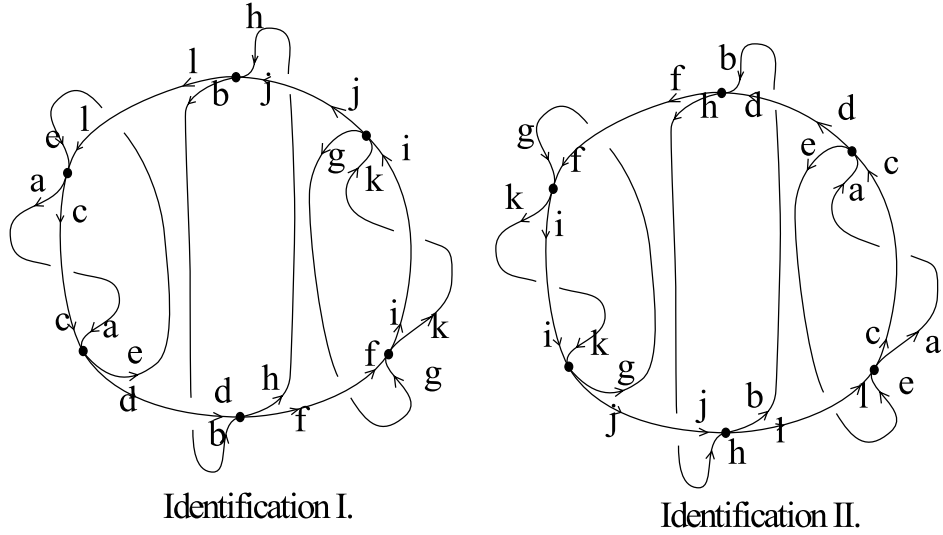}}
\label{fig:8}
\caption{ The two possible labelings. }
\end{figure}

Let us now analyse this example. First a simple Euler characteristic computation shows that the train track $\tau$ is embedded in a surface of genus 3 with 2 boundary components, that are the boundary components of the regular neighborhood of $\tau$. With the labeling of Figure 8 these two boundary components are written as cyclic edge paths on $\tau$ given by :
\centerline{
$ \partial_{1} = i j \overline{h}\overline{d} e a \overline{c}\overline{l}b f
\overline{g}\overline{k}$  and
$ \partial_{2} = i  \overline{k}\overline{g} j b \overline{d}\overline{c} a e
\overline{l}\overline{h} f$.}
The notation $ \overline{x}$ means that the oriented edge $x$ is crossed by the path with the opposite orientation.

In order to make the map $\varphi$ explicit we just have to chose one of the two possible labelings of Figure 8 and read the image of each edge as an edge path on $\tau$. Using the identification II of Figure 8 for $\tau$ and reading the paths of Figure 7 on this labelled train-track gives the following cellular map $\varphi_1$ :

 $ \bullet$  $ a \mapsto k , \hspace{2mm}$ 
$ \bullet$ $ b \mapsto f i j , \hspace{2mm}$
$ \bullet$ $ c \mapsto k g k , \hspace{2mm}$
$ \bullet$ $d \mapsto j , \hspace{2mm}$
$ \bullet$ $e \mapsto j b f , \hspace{2mm}$
$ \bullet$ $f \mapsto l a , \hspace{2mm}$\\
$ \bullet$ $g \mapsto a, \hspace{2mm}$
$ \bullet$  $h \mapsto l c d , \hspace{2mm}$
$ \bullet$ $i \mapsto e , \hspace{2mm}$
$ \bullet$  $j \mapsto  a d, \hspace{2mm}$
$ \bullet$  $k\mapsto  d h l , \hspace{2mm}$
$ \bullet$ $l \mapsto f$.

The first very good news is that no letter appears in it's own image. This is enough to prove that
 $\varphi_1$ has no fixed points by Lemma 2.4. This implies that the homeomorphism that represents 
$\varphi_1$ on the genus 3 surface with 2 boundary components has no fixed points in it's interior. We need now to check what happens for the boundaries. A simple computation shows that  $\varphi_1 ( \partial_{i}), i = 1, 2$ is a word that cyclically reduces to a cyclic permutation of  $ \partial_{i} , i = 1, 2$. Thus each boundary component is invariant under $\varphi_1$, i.e. the image is isotopic to itself, as a closed curve, but with a non trivial rotation because of the cyclic permutation.\\
This proves first that $\varphi_1$ is induced by a homeomorphism $f_1$  that has no fixed points on the genus 3 surface with 2 boundary components.

We observe that each boundary curve $\partial_{i}, i = 1, 2$, as an edge path in $\tau$, has six cusps and six sides. Checking the singularity type property using the Nielsen paths method of Proposition 2.3 is unfortunately not necessary, indeed the bad news is that the incidence matrix $M (\tau, \varphi_1)$ fails to be irreducible. Indeed, we check that the proper subgraph 
$\tau_0$ of $\tau$ that consists of the edges $\{a,c,d,f,g,h,j,k,l \}$ is invariant under 
$\varphi_1$ and therefore the matrix $M (\tau, \varphi_1)$ has a block structure and Proposition 2.2 does not apply. The existence of this invariant subgraph proves that the homeomorphism $f_1$ is in fact reducible and an invariant curve is easy to exhibit.

\subsection{Infinite sequence of pseudo-Anosov homeomorphisms.}

The idea to create a sequence of pseudo-Anosov homeomorphisms from the map $\varphi_1$ already obtained, is to compose it with an infinite sequence of Dehn twists along a curve $C$ that belongs to only one part of $\tau$, with respect to the symmetry, and do not change the map on the symmetric part of $\tau$. In other words this particular curve $C$ should have the property that it's image under the homeomorphism $f_1$ is disjoint from itself. An arbitrary train-track map will certainly not admit such a curve. Our previous example, in the labelled train-track of Figure 8  admits such a curve given by the edge path $ C = i. g$, it's image $\varphi_1 (C) = C' = e.a$ is disjoint from $C$. It turns out that the Dehn twist along $C$ is possible to describe as splitting sequence on $\tau$. With the labeling of Figure 
8 and our previous notations, such a splitting sequence is given by : 

\centerline{
$   \Big\{\frac{t (f)}{i(i)} ; \frac{i (j)}{t(i)} ;
\frac{i (k)}{t(g)}  ;  \frac{t (k)}{i(g)} \Big\}$.  }

The new train-track $\tau'$ is homeomorphic to $\tau$, with a relabeling map 
$\alpha : \tau \rightarrow \tau'$ given by:
\centerline{ $\bullet  \hspace{2mm} i \rightarrow g, \bullet \hspace{2mm} g \rightarrow i, \bullet \hspace{2mm}  x \rightarrow x$ for all $x\notin 
\{i, g\}$.}

The image of $\tau'$ is shown by Figure 9 and the twist map, along the closed curve 
$ C = i . g$,  called  $ T(i,g) : \tau ' \rightarrow \tau $ is given by:\\
\centerline{
$\bullet \hspace{2mm} f \rightarrow f.i ,\bullet  \hspace{2mm} j \rightarrow i.j ,\bullet \hspace{2mm} k\rightarrow g.k.g $  and $ \bullet  \hspace{2mm} x \rightarrow x$ for all $x \notin \{f, j, k\}$ . }

\begin{figure}[htbp]
\centerline {\includegraphics[angle=0,height=90mm]{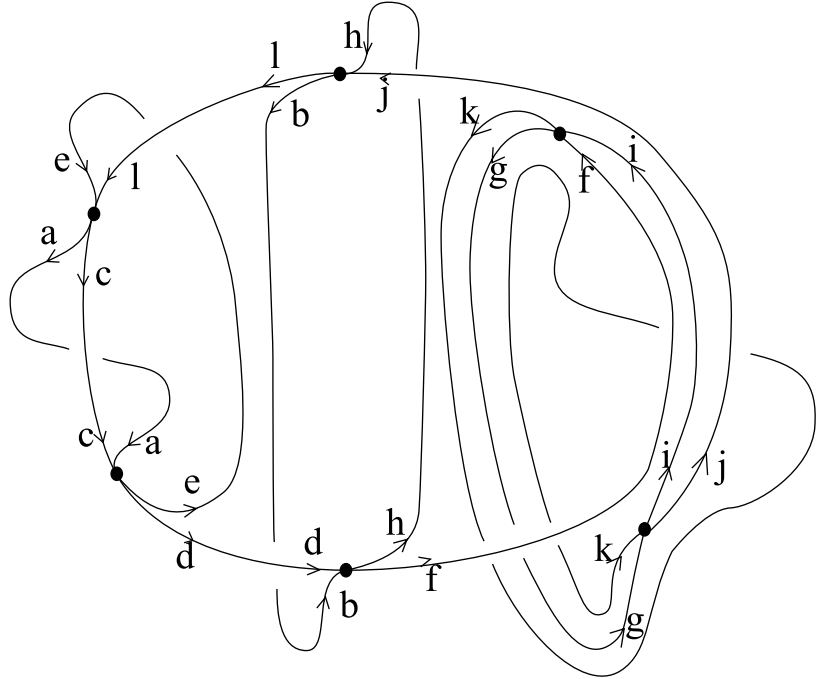} }
\label{fig:9}
\caption{  The Dehn twist $ T(i,g)$.}
\end{figure}

The composition of all these maps : 
$  \hspace{5mm} \varphi_2 : \tau '  \stackrel{ T(i,g) }{\rightarrow} \tau 
 \stackrel{ \varphi_1 }{\rightarrow} \tau 
\stackrel{\alpha}{\rightarrow} \tau ' $ gives the following : 

 $ \bullet$  $ a \mapsto k , \hspace{2mm}$ 
$ \bullet$ $ b \mapsto f g j , \hspace{2mm}$
$ \bullet$ $ c \mapsto k i k , \hspace{2mm}$
$ \bullet$ $d \mapsto j , \hspace{2mm}$
$ \bullet$ $e \mapsto j b f , \hspace{2mm}$
$ \bullet$ $f \mapsto l a e, \hspace{2mm}$\\
$ \bullet$ $g \mapsto a, \hspace{2mm}$
$ \bullet$  $h \mapsto l c d , \hspace{2mm}$
$ \bullet$ $i \mapsto e , \hspace{2mm}$
$ \bullet$  $j \mapsto e a d, \hspace{2mm}$
$ \bullet$  $k\mapsto a d h l a, \hspace{2mm}$
$ \bullet$ $l \mapsto f$.

We check that $\varphi_2$ has no fixed points by Lemma 2.4 and that the incidence matrix 
$M(\tau', \varphi_2)$ is irreducible. We now have to compute the boundary periodic points and check, via the Nielsen path method of Proposition 2.3, if the singularity type condition is satisfied in order to apply Proposition 2.2.\\
For the labelled train-track $\tau'$ (see figure 9) the two boundary components are given by:\\

\centerline{
$ \partial_{1} = c d  \overline{b}  \overline{j}  i  k   \overline{g}\overline{f}  h  l 
\overline{e}  \overline{a}$  and
$ \partial_{2} = l  c   \overline{a}  \overline{e}  d  h  \overline{j}  \overline{g}  k   i 
\overline{f}  \overline{b} $.  }

Each of these components have six cusps and six sides. For concreteness let us focus our computation on one of these sides, for instance the one denoted $[c . d]$ along $ \partial_{1}$.
We consider the successive images of this side under $\varphi_2$ in order to find the boundary periodic points and their location along  $ \partial_{1}$.\\
We first check that each of these sides have period 3 :\\

\centerline{
$[ \check{c}. d]  \rightarrow   k . [i.   \check{k}]. j  \rightarrow .... e.a.d. [ \check{h} . l] a  ....  \rightarrow .... l.   [  \check{c}. d] . f ...$ }

Where the bracket in the previous writing represents the position of one side along the boundary curve
 $ \partial_{1}$ and the symbol $ \check{..}$  shows which edge on the side contains the periodic point (here of period 3).
It is then an obvious checking that the side $[c.d]$ has only one periodic boundary point.
By a similar computation we show that this single boundary periodic point property is also satisfied for the other orbit along the sides of $ \partial_{1}$ (for instance the side 
$[ \overline{b} .  \overline{j}] $ ), as well as for the orbits of the sides along the component 
$ \partial_{2}$.

Therefore the measured foliation that is constructed from $\varphi_2$ satisfies the singularity type condition by Proposition 2.3 and  Proposition 2.2 implies that $\varphi_2$ represents a pseudo-Anosov homeomorphism either on a surface of genus 3 with 2 punctures, without fixed points, or on a closed surface of genus 3 with two fixed 6-prong singularities of positive index (the six separatrices are permuted in two orbits of period 3 from our previous computation).
In both cases the pseudo-Anosov homeomorphisms have orientable invariant foliations by Lemma 2.6. 
We thus obtain our first example.

The construction of an infinite sequence is now clear from this example, we just have to iterate the twist map. In fact since the two train-tracks $\tau$ and $\tau'$ differs just by an order two relabeling (denoted by $\alpha$ above) there is another sequence of splitting morphisms for the next twist that we denote $T(g,i) : \tau \rightarrow \tau'$ and is given by the following splitting sequence:\\

\centerline{
$\Big\{\frac{t (f)}{i(g)} ; \frac{i (j)}{t(g)} ;
\frac{i (k)}{t(i)}  ;  \frac{t (k)}{i(i)} \Big\}$.  }

The next map in the sequence is obtained by the composition:
$  \hspace{5mm} \varphi_3 : \tau  \stackrel{ T(g,i) }{\rightarrow} \tau '  \stackrel{ T(i,g) }{\rightarrow} \tau 
 \stackrel{ \varphi_1 }{\rightarrow} \tau  $.\\
Observe that the relabeling $\alpha$ is not needed here and this map $\varphi_3$ is written combinatorially on the labelled train-track of Figure 8 with identification II as :

 $ \bullet$  $ a \mapsto k , \hspace{2mm}$ 
$ \bullet$ $ b \mapsto f i j , \hspace{2mm}$
$ \bullet$ $ c \mapsto k g k , \hspace{2mm}$
$ \bullet$ $d \mapsto j , \hspace{2mm}$
$ \bullet$ $e \mapsto j b f , \hspace{2mm}$
$ \bullet$ $f \mapsto l a e a, \hspace{2mm}$\\
$ \bullet$ $g \mapsto a, \hspace{2mm}$
$ \bullet$  $h \mapsto l c d , \hspace{2mm}$
$ \bullet$ $i \mapsto e , \hspace{2mm}$
$ \bullet$  $j \mapsto a e a d, \hspace{2mm}$
$ \bullet$  $k\mapsto e a d h l a e, \hspace{2mm}$
$ \bullet$ $l \mapsto f$.

This new map $\varphi_3$ satisfies exactly the same properties than $\varphi_2$. 

The simplest infinite sequence can then be written as:
$  \hspace{5mm} \varphi_{2n+1} : \Big\{ \tau  \stackrel{ T(g,i) }{\rightarrow} \tau '  \stackrel{ T(i,g) }{\rightarrow} \tau {\Big\}}^n
 \stackrel{ \varphi_1 }{\rightarrow} \tau  $,\\
 where the notation $\{ ...\}^n$ is the natural $n$ times composition of splitting morphisms.
The combinatorial map is given on the train-track $\tau$ of Figure 8, with identification II, as:\\

 $ \bullet$  $ a \mapsto k , \hspace{2mm}$ 
$ \bullet$ $ b \mapsto f i j , \hspace{2mm}$
$ \bullet$ $ c \mapsto k g k , \hspace{2mm}$
$ \bullet$ $d \mapsto j , \hspace{2mm}$
$ \bullet$ $e \mapsto j b f , \hspace{2mm}$
$ \bullet$ $f \mapsto l a (e a)^{n}, \hspace{2mm}$\\
$ \bullet$ $g \mapsto a, \hspace{2mm}$
$ \bullet$  $h \mapsto l c d , \hspace{2mm}$
$ \bullet$ $i \mapsto e , \hspace{2mm}$
$ \bullet$  $j \mapsto (a e)^{n} a d, \hspace{2mm}$
$ \bullet$  $k\mapsto (e a)^{n} d h l (a e)^{n}, \hspace{2mm}$
$ \bullet$ $l \mapsto f$.

\begin{proposition}
The sequence $\{  \varphi_{2n+1} : \tau \rightarrow \tau ; n \in \mathbb{N}^*  \}  $ is a sequence of train-track representatives of a sequence of homeomorphisms $ f_{2n+1}: S \rightarrow S  $ such that :\\
- each $f_i$ is pseudo-Anosov on the surface $S$ of genus 3.\\
- each $f_i$ has an orientable invariant measured foliation with 2 fixed 6-prong singularities of positive index and no other fixed points.
\end{proposition}

{\em Proof.} The fixed point property is directly checked by Lemma 2.4. The pseudo-Anosov property is checked inductively with Propositions 2.2 and 2.3. The irreducibility property of the incidence matrix $M(\varphi_3, \tau)$ is checked as for $M(\varphi_2, \tau)$. Then each entry of 
$M(\varphi_{2n+1}, \tau)$ is increasing with respect $n$ and thus if $M(\varphi_3, \tau)$ is irreducible then all $M(\varphi_{2n+1}, \tau)$ are irreducible. 

The single boundary periodic point property of Proposition 2.3 for $\varphi_{2n +1}$ is checked directly  as for  $\varphi_3$. This implies first that the $ f_{2n+1}$ are all pseudo-Anosov. The invariant foliations of these pseudo-Anosov homeomorphisms are all orientable by Lemma 2.6 and are all in the same stratum because of the same 6-cusps structure of the boundary curves $\partial_1 , \partial_2$ and the single boundary periodic point property of Proposition 2.3. 
Observe that the sequence of growth rates of these pseudo-Anosov homeomorphisms goes to infinity, this is obvious, but this sequence is computable since the sequence of incidence matrices is explicit.
$\square$

For completeness, the maps $\varphi_{2n +1}$ is written, as a splitting sequence $ \mathcal{S}_{2n +1} $ with the above notations as:

$\Big\{   \big\{ \frac{i (b)}{t(l)} ; \frac{t (b)}{i(d)} ; \frac{t (k)}{i(f)} ; \frac{i (k)}{t(j)} ;
\frac{t (e)}{i(l)} ; \frac{t (c)}{i(a)} ; \frac{i (c)}{t(a)} ; \frac{i (e)}{t(d)} ;
\frac{i (h)}{t(f)} ;  \frac{t (h)}{i(j)} ; \frac{t (f)}{i(g)} ; \frac{i (j)}{t(g)} \big\}
\big\{ \frac{t (f)}{i(i)} ; \frac{i (j)}{t(i)} ;
\frac{i (k)}{t(g)}  ;  \frac{t (k)}{i(g)} .
\frac{t (f)}{i(g)} ; \frac{i (j)}{t(g)} ;
\frac{i (k)}{t(i)}  ;  \frac{t (k)}{i(i)}\big\}^{n} \Big\}$.

This completes the proof of the main theorem.$\square$

{\em Remark.} With the same basic ingredients there is another infinite sequence that is given by the composition:
$  \hspace{5mm} \psi_{n} :  \tau  \stackrel{ \varphi_1 }{\rightarrow}
\Big\{ \tau  \stackrel{ T(g,i) }{\rightarrow} \tau '  \stackrel{ T(i,g) }{\rightarrow} \tau {\Big\}}^n  $.\\
The corresponding combinatorial maps on the labelled train track of Figure 8 are given by:

$ \bullet$  $ a \mapsto (i g) ^n k (g i )^n , \hspace{2mm} $
$ \bullet$   $ b \mapsto f(i g )^n  i (g i)^n  j , \hspace{2mm}$
$ \bullet$   $ c \mapsto (i g )^n k (gi) ^n  g  (i g)^n  k  (g i) ^n,\hspace{2mm} $
$ \bullet$    $d \mapsto (g i) ^n  j , \hspace{2mm}$\\
$ \bullet$     $e \mapsto  (g i) ^n j b f (ig )^n , \hspace{2mm}$
 $ \bullet$    $f \mapsto l a , \hspace{2mm}$
$ \bullet$     $g \mapsto a, \hspace{2mm}$
$ \bullet$      $h \mapsto l c d , \hspace{2mm}$
$ \bullet$      $i \mapsto e ,\hspace{2mm}$
$ \bullet$       $j \mapsto  a d, \hspace{2mm}$
$ \bullet$     $k\mapsto  d h l , \hspace{2mm}$\\
$ \bullet$     $l \mapsto f (i g) ^n \hspace{2mm}$.\\
This sequence is different from the previous one and all these maps satisfy again the same properties with respect to the fixed points question.

 These sequences are very special in particular because the Dehn twists are directly realized as splitting sequences on the same train-track. Similar examples on other surfaces and other strata should be obtained the same way but each surface and each stratum requires a specific study. For genus two surfaces there is a work in progress in this direction.  A systematic study of all the pseudo-Anosov homeomorphisms on a given stratum would require to develop the specific train-track automata theory. For practical purposes such automata are huge and an algorithmic description is probably necessary, in particular for counting problems.

\pagebreak

{\Large{\bf References}}

[ALM]  L.Alseda, J.Llibre, M.Misiurewicz. {\em Combinatorial dynamics and entropy in dimension one.} World scientific, Advances series in non linear dynamics, Vol 5, (1993).

[A] P.Arnoux. {\em Un exemple de semi-conjugaison entre un \'echange d'intervalles et une translation du tore.} Bull. Soc. Math. France, 116 (1988), n¡4, pp 489-500.

Ê[A.Y] P. Arnoux , J.-C.Yoccoz. {\em ÊConstruction de diffŽomorphismes pseudo-Anosov.}Ê
C. R. Acad. Sci. Paris SŽr. I Math.Ê292Ê(1981), no. 1, 75--78.Ê

[Ba] M. Bauer. {\em Examples of pseudo-Anosov homeomorphisms}, Trans. A.M.S, {\bf 330}
1992, no. 1, pp 333-359.

[BH] M.Bestvina, M.Handel. {\em Train tracks for surface homeomorphisms.}
Topology, {\bf 34} (1995) 109-140.

[Bo] P.Boyland. {\em Topological methods in surface dynamics.} Top. Appl.  58, ( 1994), pp 223-298.

[Bu] A.Bufetov. {\em Logarithmic asymptotics for the number of periodic orbits of the Teichm\"uller flow on Veech's space of zippered rectangles.} Moscow.Math.Journal, To appear.

[CB] A.Casson, S.Bleiler. {\em Automorphisms of surfaces after Nielsen and Thurston} Cambridge University Press, Cambridge (1988).\\
Handwritten notes, University of Texas at Austin, volume 1 Fall 1982, Volume 2, Spring 1983.

[EM] A.Eskin, M.Mirzakhani. {\em Counting closed orbits in Modular space.} Preprint.

[FLP] A.Fathi, F.Laudenbach, V.Poenaru. {\em Travaux de Thurston sur les surfaces} Ast\'erisque vol 66-67 (1979).

[FL] J.Fehrenbach, J.Los. {\em Roots, symmetry and conjugacy of pseudo-Anosov homeomorphims.} Preprint 2007.

[Ha] U.Hamenst\"adt. {\em Train tracks and the Gromov boundary of the complex of curves} in "Spaces of Kleinian groups"  (Y.Minsky, M.Sakuma, C.Series, Eds.)
Lond. Math.Soc..Lec.Notes 329 (2005).

[J] B.Jiang. {Lectures on Nielsen fixed point theory.} A.M.S.  Providence RI (1983).

[KLS] K.H.Ko, J.Los, W.T.Song. {\em Entropy of braids.} Jour. Knot. Theo. ramif,
11, 4, (2002), 647-666.

[Lo] J.Los. {\em On the forcing relation for surface homeomorphisms.} Publ.Math. IHES, 85, (1997), 5-61.

[Mo86] L. Mosher. {\em The classification of pseudo-Anosov's , } Low-dimensional topology and Kleinian groups (Coventry/Durham, 1984), LMS Lecture Notes Series, Vol 112, (1986), pp 13-75.

[Mo03] L.Mosher. {\em Train-track expansions of measured foliations.} Preprint monograph 2003.

[PaPe] A.Papadopoulos,  R.Penner. {\em A characterization of Pseudo-Anosov foliations.}  Pac. Journ.Math. 130, 2 (1987), pp 359-377.

[PH] R.C.Penner with J.L.Harer. {\em Combinatorics of train tracks}
Annals of Math.Studies, 125, (1992).

[Ra] G.Rauzy. {\em Echanges d'intervalles et transformations induites.}
Acta.Arith. 34, (1979), 315-328.

[Th] W. Thurston. {\em On the geometry and dynamics of diffeomorphisms on surfaces}, Bull. A.M.S. {\bf 19} (1988) pp. 417-431.

[Ve]  W.A.Veech. {\em Teichm\"uller geodesic flow.} Annals of Math. 124, (1986), 441-530.

[Wi] R.Williams. {\em One dimensional non wandering sets.}
 Topology, 6, ( 1967), 473-487.

\bigskip

Laboratoire d'Analyse, Topologie, Probabilit\'e (LATP), UMR CNRS 6632\\
Universit\'e de Provence,\\
39 rue F.Joliot Curie\\
13453 Marseille Cedex 13,\\
France.\\
E-mail adresse : los@cmi.univ-mrs.fr

\end{document}